\documentclass[12pt]{article}
\usepackage{amssymb,amsmath}
\def\hang{\hangindent\parindent}
\def\tex#1{\indent\llap{[#1]\enspace}\ignorespaces}
\def\re{\par\hang\tex}
\def\F{{\mathbb{F}}}
\def\a{\alpha}
\def\b{\beta}
\def\d{\delta}

\def\BB{\mathcal B}

\def\l{\lambda}
\def\L{\Lambda}

\def\sc{\scriptstyle}
\def\ssc{\scriptscriptstyle}

\def\cl{\centerline}

\def\ul{\underline}

\def\vs{\vspace*}

\def\AA{{\mathcal A}}

\def\PP{{\mathcal P}}

\def\N{\mathbb{N}{\ssc\,}}
\def\Z{\mathbb{Z}{\ssc\,}}

\numberwithin{equation}{section}
\newtheorem{theo}{Theorem}[section]

\newtheorem{lemm}[theo]{Lemma}

\newtheorem{clai}{Claim}

\begin{document}
\cl {\Large\bf Highest weight representations of a} \cl {{\Large\bf
Lie algebra of Block type}\footnote {Supported by a NSF grant
10471096 of China, ``One Hundred Talents Program'' from University
of Science and Technology of China.}} \vs{6pt} \cl{Yuezhu Wu
$^{1,2)}$,  Yucai Su $^{3)}$} \cl{\small $^{1)}$Department of
Mathematics, Shanghai Jiao Tong University,
 Shanghai 200240, China}

\cl{\small $^{2)}$Department of Mathematics, Qufu Normal
University,
 Qufu 273165, China}

 \cl{\small \small $^{3)}$Department of
Mathematics, University of Science and Technology of
\vs{-3pt}China} \cl{\small Hefei 230026, China}

 \cl{\small E-mail:  ycsu@ustc.edu.cn} \vs{6pt}
\vs{10pt}\par

{\small
\parskip .005 truein
\baselineskip 3pt \lineskip 3pt

\noindent{{\bf Abstract.} For a field $\F$ of characteristic zero
and an additive subgroup $G$ of $\F$, a Lie algebra $\BB(G)$ of
Block type is defined with basis $\{L_{\a,i},\,c\,|\,\a \in
G,\,-1\le i\in\Z\}$ and relations
$[L_{\a,i},L_{\b,j}]=((i+1)\b-(j+1)\a)L_{\a+\b,i+j}+\a\d_{\a,-\b}\d_{i+j,-2}c,\,
[c,L_{\a,i}]=0.$ Given a total order $\succ$ on $G$ compatible with
its group structure, and any $\L\in \BB(G)_0^*$, a Verma
$\BB(G)$-module $M(\L,\succ)$ is defined, and the irreducibility of
$M(\L,\succ)$ is completely determined. Furthermore, it is proved
that an irreducible highest weight $\BB(\Z)$-module is quasifinite
if and only if it is a proper quotient of a Verma module. \vs{5pt}

\noindent{\bf Key words:} Verma modules, Lie algebras of Block
type, irreducbility.}

\noindent{\it Mathematics Subject Classification (2000):} 17B10,
17B65,  17B68.}
\parskip .001 truein\baselineskip 8pt \lineskip 8pt

\vs{6pt}

\par

\cl{\bf\S1. \
Introduction}\setcounter{section}{1}\setcounter{equation}{0} Block
[B] introduced a class of infinite dimensional simple Lie algebras
over a field of characteristic zero. Generalizations of Block
algebras (usually referred to as {\it Lie algebras of Block type})
have been studied by many authors (see, for example, [DZ, LT, S1,
S2, X1, X2, WZ, ZM]). Partially because they are closely related to
the Virasoro algebra (and some of them are sometimes called
Virasoro-like algebras), these algebras have attracted some
attention in the literature.

Let $\F$ be a field of characteristic 0 and $G$  an additive
subgroup of $\F$.
 The {\it Lie algebra $\BB(G)$ of Block type} considered in this paper is  the Lie algebra with basis
$\{L_{\a,i},c\,|\,\a \in G,\, i\in\Z$, $i\geq -1\}$ and relations
\begin{equation}
[L_{\a,i},L_{\b,j}]  =
((i+1)\b-(j+1)\a)L_{\a+\b,i+j}+\a\d_{\a,-\b}\d_{i+j,-2}c,
~~[c,L_{\a,i}] = 0.\label{def}
\end{equation}
Let
\begin{equation}\BB(G)_\a={\rm span}\{L_{\a,i}|~i\geq
-1\}+\d_{\a,0}\F c. \label{grad}
\end{equation}
 Then $\BB(G)=\oplus_{\a\in
G}\BB(G)_\a$ is $G$-graded (but not finitely graded). Throughout
this paper, we fix a total order `` $\succ $'' on $G$ compatible
with its group structure. Denote
$$ G_{+}=\bigl\{x\in G\bigm|x\succ 0\bigr\},\quad
G_{-}=\bigl\{x\in G\bigm|x\prec 0\bigr\}. $$ Then
$G=G_{+}\cup\bigl\{0\bigr\}\cup G_{-}$. Setting
$\BB(G)_{\pm}=\oplus_{\pm\a\succ 0}\BB(G)_\a,$ we have the
triangular decomposition
$$\BB(G)=\BB(G)_-\oplus \BB(G)_0 \oplus \BB(G)_+.$$
Note that $\BB(G)_0={\rm span}\{L_{0,i}\,|~i\geq -1\}\oplus\F c$ is
a commutative subalgebra of $\BB(G)$ (but it is not a Cartan
subalgebra).

A $\BB(G)$-module $V$ is {\it quasifinite} if $V$ is finitely
$G$-graded, namely, $$V =\oplus_{\a\in G}V_\a,\ \ \ \BB(G)_\a
V_\b\subset V_{\a+\b},\ \ {\rm dim}V_\a<\infty \mbox{\ \ \ for \ \
}\a,\b\in G.$$ Quasifinite modules are closed studied by some
authors, e.g., [KL, KR, S1, S2]. In [S1], it is proved that a
quasifinite irreducible $\BB(\Z)$-module is a highest or lowest
weight module and the quasifinite irreducible highest weight modules
are classified. The main result of this paper is the following.
\begin{theo}
\label{theo1}\vskip-6pt
\begin{itemize}\parskip-6pt
\item[{\rm(1)}]
An irreducible highest weight $\BB(\Z)$-module is quasifinite if and
only if it is a proper quotient of a Verma module.
\item[{\rm(2)}]
Let $\L\in \BB(G)_0^*$. With respect to a dense order $`` \succ$" of
$G$ $($cf.~$(\ref{dense}))$, the Verma $\BB(G)$-module $M(\L,\succ)$
is irreducible if and only if $\L\neq 0$. Moreover, in case $\L=0$,
if we set
$$M'(0,\succ)=\sum\limits_{k>0,\, \a_1,\cdots,\a_k\in G_+}\F
L_{-\a_1,i_1}\cdots L_{\a_k,i_k}v_0,$$ then $M'(0,\succ)$ is an
irreducible submodule of $M(0,\succ)$ if and only if for all $x,y\in
G_+,$ there exists a positive integer $n$ such that $nx\succ y.$

 \item[{\rm(3)}] With respect to a
discrete order $`` \succ$" $($cf.~$(\ref{discrete}))$, the Verma
$\BB(G)$-module $M(\L,\succ)$ is irreducible if and only if
$M_a(\L,\succ)$ is an irreducible $\BB(a\Z)$-module.
\end{itemize}
\end{theo}

\cl{\bf\S2. \ Verma modules over $\BB(G)$
}\setcounter{section}{2}\setcounter{theo}{0}\setcounter{equation}{0}Let
$U=U(\BB(G))$ be the universal enveloping algebra of $\BB(G)$. For
any $\L\in \BB(G)_0^*$ (the dual space of $\BB(G)_0$), let
$I(\L,\succ)$ be the left ideal of $U$ generated by the elements
$$\{L_{\a,i}\,|~\a\succ 0, i\geq -1\}\cup\{h-\L(h)\cdot 1\,|~h\in \BB(G)_0\}.$$
Then the {\it Verma $\BB(G)$-module} with respect to the order
``$\succ$''
  is defined as $$M(\L,\succ)=U/I(\L,\succ).$$
By the PBW theorem, it has  a basis consisting of all vectors of the
form
$$L_{-\a_1,i_1}L_{-\a_2,i_2}\cdots L_{-\a_k,i_k}v_{\L},$$
where $v_\L$ is the coset of $1$ in $M(\L,\succ)$, and
$$-1\le i_j\in \Z,\ \ \ 0\prec\a_1\preccurlyeq
\cdots\preccurlyeq \a_k, \mbox{ \ \ and \ }i_s\leq i_{s+1}\mbox{ \
if \ }\a_s=\a_{s+1}.$$ Note that $M(\L,\succ)$ is a highest weight
$\BB(G)$-module  in the sense that
$M(\L,\succ)=\oplus_{\a\preccurlyeq 0}M_{\a},$ where $M_0=\F v_\L,$
and $M_\a$ for $\a\prec 0$ is spanned by
\begin{equation}\label{weight}
L_{-\a_1,i_1}L_{-\a_2,i_2}\cdots L_{-\a_k,i_k}v_\L,\end{equation}
with $i_j\geq -1,\, 0\prec\a_1\preccurlyeq \cdots \preccurlyeq
\a_k,$ and $\a_1+\cdots+ \a_k=-\a.$ Thus $M(\L,\succ)$ is a
$G$-graded $\BB(G)$-module with ${\rm dim}M_{-\a}=\infty$ for any
$\a \in G_+.$

 We call a nonzero vector $u\in M_\a$  a {\it weight vector with weight $\a$.}
 For any $a\in G$,
denote
$$\BB(a\Z)={\rm span}\{L_{na,k}\,|~a\in \Z,k\geq -1\},$$
a subalgebra of $\BB(G)$ isomorphic to $\BB(\Z)$. We also denote
$$M_a(\L,\succ)\ \ =\ \ \mbox{ $\BB(a\Z)$-submodule of $M(\L,\succ)$
 generated by $v_\L$.}$$
Denote
\begin{equation}\label{denote-Ba}
B(\a)=\{\b\in G\,|~0\prec \b \prec \a\}\mbox{ \ \
 for \ }\a\in G_+.\end{equation}
 The order `` $\succ $" is called {\it dense} if
\begin{equation}\label{dense}
\# B(\a)=\infty \mbox{ \ \ for all \  }\a \in G_+,\end{equation} it
is  {\it discrete} if
 \begin{equation}\label{discrete}
 B(a)=\emptyset\mbox{ \ \ for some \ }a\in
 G_+.\end{equation}
\vskip5pt \noindent{\it Proof of Theorem \ref{theo1}(2) and
(3).~}~(2) Suppose the order $`` \succ$" is dense. For each
$m\in\N=\{1,2,...\}$, set
\begin{equation}
V_{m}=\sum\limits_{\stackrel{^{^{\sc 0\leq s\leq m,\ i_1,\cdots
i_s\geq -1}}}{_{_{\sc 0\prec\a_{1}\preccurlyeq\cdots\preccurlyeq
\a_s}}}} \F L_{-\a_1,i_1}\cdots L_{-\a_s,i_s}v_{\L},
\label{gradation}
\end{equation} where $i_s\leq i_{s+1}\mbox{ if
}\a_s=\a_{s+1}.$ It is clear that $L_{\a,k}V_m\subseteq V_m$ for any
$\a\in G_{+},\, k\geq -1$.

Let $u_0\neq 0$ be any given weight vector in $M(\L,\succ)$. We want
to prove that $v_{\L}\in U(\BB(G))u_{0}$ if $\L\ne0$. We divide the
proof into four steps:
\smallskip

{\it Step 1.\,} We claim that there exists some weight vector  $u\in
U(\BB(G))u_0$ 
such that, for some
$r\in\Bbb{N}$,
$$
u\equiv \sum\limits_{k_1,\cdots,k_r\in
\N}a_{\ul{k}}L_{-\varepsilon_r,k_r}\cdots L_{-\varepsilon_1,k_1}v_\L
\mbox{ } ({\rm mod}\mbox{ } V_{r-1})\mbox{ \ \ for some \ }a_{\ul
k}\in\F,
$$ where $0\prec\varepsilon_r\prec\cdots\prec\varepsilon_1$, and $ 0\ne a_{\ul{k}}\in\F$
 for some $\ul{k}=(k_r,\cdots,k_1).$

It is clear that $u_0\in V_{r}\setminus{V_{r-1}}$ for some
$r\in\Bbb{N}$. If $r\le 1$, our claim clearly holds. So we assume
that $r>1$.  Hence we can write
$$ u_0\equiv \sum\limits_{0\prec\a_{1}\preccurlyeq\cdots\preccurlyeq \a_r }
a_{\ul{\a},\ul{i}}L_{-\a_1,i_1}\cdots L_{-\a_r,i_r}v_\L\mbox{ }
({\rm mod}\mbox{ } V_{r-1})\mbox{ \ \ for some \
}a_{\ul{\a},\ul{i}}\in\F,
$$ where $\ul\a=(\a_1,...,\a_r),\,\ul i=(i_1,...,i_r)$, and we
denote $$(\ul{\a},\ul{i})=(\a_1,...,\a_r,i_1,...,i_r).$$ Let
$I=\bigl\{(\ul{\a},\ul{i})\,|\,a_{\ul{\a},\ul{i}}\neq 0\bigr\}$. By
assumption, $I\neq\emptyset$. For any $\ul{\a}$ and
$\ul{\a}'=(\a'_1,\cdots,\a'_r),$ we define
\begin{equation}\label{order1} \ul{\a}\succ \ul{\a}'\,\ \ ~ \Longleftrightarrow \ \ \,
\exists\,s\in\{1,...,r\}\mbox{ \ such that \ } \a_s\succ\a'_s,\mbox{
\ and \ } \a_t=\a'_t\mbox{ \ for \ }t>s .
\end{equation}
Similarly, for any $\ul{i}$ and $\ul{i}'=(i'_1,\cdots,i'_r),$ we
define
\begin{equation}\label{order-i} \ul{i}> \ul{i}'
\,\ \ ~ \Longleftrightarrow \ \ \, \exists\,s\in\{1,...,r\}\mbox{ \
such that \ } i_s>i_s,\mbox{ \ and \ } i_t=i'_t\mbox{ \ for \ }t>s .
\end{equation}
  For any
$(\ul{\a},\ul{i}),\,(\ul{\a}',\ul{i}')\in I$, we define
\begin{equation}\label{order1+} (\ul{\a},\ul{i})\succ(\ul{\a}',\ul{i}')\, \ \ \
\Longleftrightarrow\,\ \ \ \ul{\a}\succ \ul{\a}'\mbox{, \ \ or \ \
}\ul{\a}=\ul{\a}',\ \ul{i}>\ul{i}'.
\end{equation}
Let $$ (\ul{\b},\ul{j})=(\b_1,\cdots,\b_r,j_1,\cdots,j_{r}),\ \
\,0\prec\b_{1}\preccurlyeq\cdots\preccurlyeq \b_{r},
$$
be the unique maximal element in $I$. Then we can write $\b$ as $$
\b=(\b_1,...,\b_s,\b_r,...,\b_r)\mbox{ \ \ for some \
}s\in\{1,...,r\}.$$ By the assumption that $`` \succ"$ is a dense
order, we can always find some $\varepsilon_1\in G_{+}$ such that
$$\varepsilon_{1}\prec \b_1\mbox{ \ and \ }\{x\in G_+\,|\,\b_{r}-\varepsilon_{1}\prec x\prec
\b_{r}\}\cap\{\a_{r}, \a_{r-1}\,|\,(\ul{\a},\ul{i})\in I \mbox{ for
some }\ul{i}{\ssc\,}\}=\emptyset.
$$
Using relations (\ref{def}), and noting that
$\b_r-\varepsilon_1-\a_k\succ 0$ if $\a_k\neq \b_r,
k\in\{1,\cdots,r\}$, by choosing $k_1\gg 0$ with $k_1\in \N,$  we
see that

$$u_{1}:=L_{\b_{r}-\varepsilon_1,k_1}u_{0}\equiv\sum\limits_{0\prec\a_{1}\preccurlyeq\cdots\preccurlyeq \a_{r-1},\,
k_1'\in\N}
a^{(1)}_{\ul{\a},\ul{j}}L_{-\varepsilon_1,k'_1}L_{-\a_1,i_1}\cdots
L_{-\a_{r-1},i_{r-1}}v_{\L}\mbox{ } ({\rm mod}\mbox{ } V_{r-1}),$$
for some $a^{(1)}_{\ul{\a},\ul{j}}\in\F$.
 Set $$
I^{(1)}=\bigl\{\bigl(\varepsilon_1,\a_1,\cdots,\a_{r-1},k'_1,i_{1},\cdots,i_{r-1}
\bigr)\bigm|a_{\ul{\a},\ul{i}}^{(1)}\neq
0\bigr\} .$$
The coefficient corresponding to
$$(\ul{\b}^{(1)},\ul{j}^{(1)})=(\varepsilon_1,\b_1,\cdots \b_{r-1},k_1+j_{s+1},j_1,\cdots,j_{r-1})$$
($(\ul{\b}^{(1)},\ul{j}^{(1)})$ maybe not the  maximal element in
$I^{(1)}$) is
$$-m((k_1+1)\b_r+(j_{s+1}+1)(\b_r-\varepsilon_1))a_{\ul{\b},\ul{j}}\neq
0\mbox{ \ \ for some \ }m\in\N.$$ Thus $ I^{(1)}\neq \emptyset.$

Now for $p=2,\cdots,r$, we define recursively and easily prove by
induction  that

(i) Let $\varepsilon_p\in G_{+}$ such that
$\varepsilon_{p}\prec\varepsilon_{p-1}$ and $$ \bigl\{x\in
G\,|\,\b_{r-p+1}-\varepsilon_{p}\prec x\prec
\b_{r-p+1}\bigr\}\cap\bigl\{\a_{r-p+1}, \a_{r-p}\,|\,
(\ul{\a},\ul{j})\in I^{(p-1)}\bigr\}=\emptyset.
$$

(ii) Choose $k_p\gg 0$ and let
$u_{p}=L_{\b_{r-p+1}-\varepsilon_p,k_p}u_{p-1}$. Then, for some
$a_{\ul{\a},\ul{j}}^{(p)}\in\F$,
$$
u_{p}\equiv\sum\limits_{{0\prec\a_{1}\preccurlyeq\cdots\preccurlyeq
\a_{r-p}}} a_{\ul{\a},\ul{j}}^{(p)}L_{-\varepsilon_{p},k'_p}\cdots
L_{-\varepsilon_1,k'_1}L_{-\a_1,i_1}\cdots
L_{-\a_{r-p},i_{r-p}}v_\L\mbox{ } ({\rm mod}\mbox{ } V_{r-1}).
$$

(iii) Let
$$I^{(p)}=\bigl\{(\varepsilon_{p},\cdots,\varepsilon_1,\a_1,\cdots, \a_{r-p},k'_p,\cdots,k'_1,i_{1},\cdots,
j_{r-p})\,|\,a_{\ul{\a},\ul{j}}^{(p)}\neq 0\bigr\}.$$ Then
$I^{(p)}\neq\emptyset$.

Now our claim follows immediately by letting $p=r$.

\smallskip
{\it Step 2.\,} We claim that there exists some weight vector  $u\in
U(\BB(G))u_0$ such that, for some $r\in\Bbb{N}$,
$$
u\equiv L_{-\varepsilon_r,k_r}\cdots L_{-\varepsilon_1,k_1}v_\L
\mbox{ } ({\rm mod}\mbox{ } V_{r-1}), $$ where $k_j\geq -1,$ and
$0\prec\varepsilon_r\prec\cdots\prec\varepsilon_1.$

By Step 1, there exists some weight vector $u\in U(\BB(G))u_0$ such
that, for some $r\in\Bbb{N}$,
$$
u\equiv \sum\limits_{k_1,\cdots,k_r\in
\N}a_{\ul{k}}L_{-\varepsilon_r,k_r}\cdots L_{-\varepsilon_1,k_1}v_\L
\mbox{ } ({\rm mod}\mbox{ } V_{r-1})\mbox{ \ \ for some \
}a_{\ul{k}}\in\F,
$$ where  $0\prec\varepsilon_r\prec\cdots\prec\varepsilon_1$ and
$$K:=\{\ul{k}=(k_r,\cdots,k_1)\,|\, a_{\ul{k}}\neq 0\}\neq
\emptyset.$$
 Let $$\ul{j}=(-1,\cdots,-1,j_s,\cdots,j_1)\mbox{ \ \  with \ }j_s\neq -1,$$
 be the unique maximal
 element in $K$ (recall (\ref{order-i})). Assume that $K$
 is not a singleton. Then $\ul{j}\neq
 (-1,\cdots,-1)$.
Set $$\d={\rm min}\{\{\varepsilon_i,
\varepsilon_j-\varepsilon_i\,|\,1\leq j<i\leq r\}\cap G_+\}.$$ Let
$\varepsilon'\in G_+$ such that $\varepsilon'\prec \d.$ Then
\begin{eqnarray*}
L_{\varepsilon',-1}\cdot u&\equiv&\sum\limits_{k_r,\cdots,k_1\in
\N}\sum\limits_{j=0}^{r}-(k_j+1)\varepsilon'a_{\ul
k}L_{-\varepsilon_r,k_r}
\cdots L_{-\varepsilon_{j+1},k_{j+1}}\\
& &\times \ L_{\varepsilon'-\varepsilon_j,k_j-1}
L_{-\varepsilon_{j-1},k_{j-1}}\cdots L_ {-\varepsilon_1,k_1}v_\L
\mbox{ } ({\rm mod}\mbox{ } V_{r-1}).
\end{eqnarray*}
The term $$L_{-\varepsilon_r,-1}\cdots L_{-\varepsilon_{s+1},-1}
L_{\varepsilon'-\varepsilon_s,j_s-1}
L_{-\varepsilon_{s-1},j_{s-1}}\cdots L_ {-\varepsilon_1,j_1}v_\L$$
appears in $L_{\varepsilon,-1}\cdot u$, since the corresponding
coefficient is $-(j_s+1)a_{\ul{j}}\neq 0.$ Using the same arguments
as above and the induction on ${\rm max}\{k_r+\cdots+
k_1\,|\,a_{\ul{k}}\neq 0\}$, we see that there exists some weight
vector  $u\in U(\BB(G))u_0$ such that
$$ u\equiv \sum\limits_{0\prec\a_{1}\preccurlyeq\cdots\preccurlyeq \a_r }
a_{\ul{\a}}L_{-\a_1,-1}\cdots L_{-\a_r,-1}v_\L\mbox{ } ({\rm
mod}\mbox{ } V_{r-1})\mbox{ \ \ for some \ }a_{\ul{\a}}\in\F.
$$ Using the same
arguments as in Step 1, we can prove the claim.

\smallskip
{\it Step 3.\,} We claim that there exists some $\varepsilon\in
G_{+}$   such that $L_{-\varepsilon,k}v_\L\in U(\BB(G))u_0$ for
$k\geq -1.$
\smallskip

By Step 2, there is a weight vector  $u\in U(\BB(G))u_0$ such that
$$ u=L_{-\varepsilon_r,k_r}\cdots
L_{-\varepsilon_1,k_1}v_\L+\sum\limits_{{\sc 0\leq
l<r,\,}{0\prec\a_{1} \preccurlyeq\cdots\preccurlyeq \a_{l}}}
b_{\ul{\a},\ul{i}}L_{-\a_1,i_1}\cdots L_{-\a_l,i_l}v_{\L}, $$ for
some $b_{\ul{\a},\ul{i}}\in\F$, where
$0\prec\varepsilon_r\prec\cdots\prec\varepsilon_1$. Assume that
$u\notin\F v_\L$.

Set
$$I_{0}=\bigl\{\ul{\a}=(\a_1,\cdots,\a_l)\,|\,b_{\ul{\a},\ul{i}}\neq
0\mbox{ \ for some \ } \ul{i}\bigr\},\ \ \ \ul{j}^{(0)}={\rm
min}\{\varepsilon_r,\a_1\,|\, \ul{\a}\in I_{0}\}.$$ Let
$\varepsilon\in G_{+}$ such that $\varepsilon\prec\ul{j}^{(0)}$.
Assume that $u\in M_\l$ (cf.~(\ref{weight})). By relations
(\ref{def}), we have
\begin{eqnarray*}
L_{-\lambda-\varepsilon,j}u &=&
f(-\lambda-\varepsilon)L_{-\varepsilon,j+k_r+\cdots+k_1}v_{\L}+\mbox{$\sum\limits_{
1\leq l<r,\, 0\prec\a_{1} \preccurlyeq\cdots\preccurlyeq \a_{l}}$}
b_{\ul{\a},\ul{i}}g_{\ul{\a},\ul{i}}(-\lambda-\varepsilon)L_{-\varepsilon,j+i_1+\cdots +i_l}v_{\L}\\
&=&\bigl\{f(-\lambda-\varepsilon)+\mbox{$\sum\limits_{\stackrel{\sc
1\leq l<r,\, 0\prec\a_{1} \preccurlyeq\cdots\preccurlyeq
\a_{l}}{i_1+\cdots +i_l=k_1+\cdots+k_r}}$}
b_{\ul{\a},\ul{i}}g_{\ul{\a},\ul{i}}(-\lambda-\varepsilon)\bigr\}L_{-\varepsilon,j+k_r+\cdots+k_1}v_{\L}\\
&&+\mbox{$\sum\limits_{\stackrel{\sc 1\leq l<r,\, 0\prec\a_{1}
\preccurlyeq\cdots\preccurlyeq \a_{l}}{i_1+\cdots+i_l\neq
k_1+\cdots+k_r}}$}
b_{\ul{\a},\ul{i}}g_{\ul{\a},\ul{i}}(-\lambda-\varepsilon)L_{-\varepsilon,j+i_1+\cdots
+i_l}v_{\L}\ \ \in \ \ U(\BB(G))u_0,
\end{eqnarray*} where in general $f(x)$ and $g_{\ul a,\ul i}(x)$ are
determinants:
\begin{eqnarray*}
\!\!\!\!&\!\!\!\!&f(x)=
         \left|
               \begin{array}{cc}
               j{\sc\!}+{\sc\!}1 & k_r{\sc\!}+{\sc\!}1\\
               x & -\varepsilon_r
               \end{array}
         \right|\
          \left|
               \begin{array}{cc}
               j{\sc\!}+{\sc\!}k_r{\sc\!}+{\sc\!}1 & k_{r-1}{\sc\!}+{\sc\!}1\\
               x-\varepsilon_r & -\varepsilon_{r-1}
               \end{array}
         \right|
         \cdots
\left|
               \begin{array}{cc}
               j{\sc\!}+{\sc\!}k_r{\sc\!}+{\sc\!}\cdots{\sc\!} +{\sc\!}k_2{\sc\!}+{\sc\!}1 & k_1{\sc\!}+{\sc\!}1\\
               x-\varepsilon_r-\cdots -\varepsilon_2 &
               -\varepsilon_1
               \end{array}
         \right|
       ,\\
\!\!\!\!&\!\!\!\!& g_{\ul a,\ul i}(x)=
         \left|
               \begin{array}{cc}
               j{\sc\!}+{\sc\!}1 & i_1{\sc\!}+{\sc\!}1\\
               x & -\a_1
               \end{array}
         \right|\
         \left|
               \begin{array}{cc}
               j{\sc\!}+{\sc\!}i_1{\sc\!}+{\sc\!}1 & i_2{\sc\!}+{\sc\!}1\\
               x-\a_1 & -\a_2
               \end{array}
         \right|
                \cdots
\left|
               \begin{array}{cc}
               j{\sc\!}+{\sc\!}i_1{\sc\!}+{\sc\!}\cdots{\sc\!}+{\sc\!}i_{k-1}{\sc\!}+{\sc\!}1 & i_k{\sc\!}+{\sc\!}1\\
               x-\a_1-\cdots-\a_{k-1} & -\a_k
               \end{array}
         \right|
       ,
\end{eqnarray*}
 Since $\deg f(x)=r>\deg g_{\ul{\a},\ul{i}}(x)$ for all $\ul{\a}\in I_{0},$
 we can find $\varepsilon\in G_{+}$ with
$\varepsilon\prec\ul{j}^{(0)}$ such that $$
f(-\lambda-\varepsilon)+\sum\limits_{ 1\leq l<r,\,0\prec\a_{1}
\preccurlyeq\cdots\preccurlyeq \a_{l}}
b_{\ul{\a},\ul{i}}g_{\ul{\a},\ul{i}}(\lambda-\varepsilon)\neq 0.
$$ So we obtain some vector
$$u=(a_1L_{-\varepsilon,i_1}+\cdots+a_nL_{-\varepsilon,i_n})v_\L\in U(\BB(G))u_0,$$
for some $0\ne a_1,\cdots,a_n\in\F.$ Choosing $\varepsilon'\in G_+$
with $\varepsilon'\prec \varepsilon,$ using
$$L_{\varepsilon-\varepsilon',-1}u=-(\varepsilon-\varepsilon')(a_1(i_1+1)L_{-\varepsilon',i_1-1}+
\cdots +a_n(i_n+1)L_{-\varepsilon',i_n-1})v_\L\in U(\BB(G))u_0,$$
and induction on $\mbox{max}\{i_1,\cdots,i_n\}$, one can deduce that
there exists some $\varepsilon'\prec \varepsilon$ such that
$L_{-\varepsilon',-1}v_\L\in U(\BB(G))u_0.$ Let $\varepsilon''\in
G_+$ such that $\varepsilon''\prec \varepsilon'.$
 Then $$L_{-\varepsilon'',k-1}v_\L=-((k+1)\varepsilon')^{-1}
 L_{\varepsilon'-\varepsilon'',k}L_{-\varepsilon',-1}v_\L
 \in U(\BB(G))u_0\mbox{ \ \ for all \ }k\geq 0.$$
This proves our claim.

\smallskip
{\it Step 4.\,} We claim that if there exists some $\varepsilon\in
G_+$ such that $L_{-\varepsilon,k}v_\L\in U(\BB(G))u_0$ for all
$k\geq -1,$ then $L_{-x,k}v_\L\in U(\BB(G))u_0$ for all  $k\geq -1$
and all $x\in B'(\varepsilon)$, where $B'(\varepsilon)$ is defined
by
$$B'(\varepsilon)={\rm span}_{\Z_+} \{y\in G_+\,|~y\preccurlyeq \varepsilon \}. $$
\smallskip

Let $\varepsilon'\in G_+$ such that $\varepsilon'\prec
\varepsilon.$
 Then $$L_{-\varepsilon',k-1}v_\L=-((k+1)(\varepsilon-\varepsilon'))^{-1}
 L_{\varepsilon-\varepsilon',-1}L_{-\varepsilon,k}v_\L
 \in U(\BB(G))u_0.$$ Since $$(k+1)\varepsilon'L_{-(\varepsilon+\varepsilon'),k-1}v_\L
 =L_{-\varepsilon',-1}L_{-
 \varepsilon,k}v_\L-L_{-\varepsilon,k}L_{-\varepsilon',-1}v_\L\in U(\BB(G))u_0,$$
 it follows that $L_{-(\varepsilon'+\varepsilon),k}v_\L\in U(\BB(G))u_0 \mbox{  for all  } k\geq -1.$
 Similarly, we deduce that $$L_{-x,k}v_\L\in U(\BB(G))u_0 \mbox{ \ \ for all \ } k\geq -1\mbox{ \ and all \ }
 x\in \Z_+\varepsilon+\Z_+\varepsilon'.$$ Our claim follows.
 \smallskip

By Step 3, we have $L_{-\varepsilon',-1}v_\L\in U(\BB(G))u_0$
 for some $\varepsilon'\in G_+.$ From
\begin{eqnarray*}
L_{\varepsilon',-1}L_{-\varepsilon',-1}v_\L
&=&\varepsilon'c\cdot v_\L=\L(c)v_\L,\\
L_{\varepsilon',k}L_{-\varepsilon',-1}v_\L&=&-(k+1)
\varepsilon'L_{0,k-1}v_\L=-(k+1)\varepsilon'\L(L_{0,k-1})v_\L\mbox{
for }k\geq 0,
\end{eqnarray*}
it is easy to see that $v_\L\in U(\BB(G))u_0$ if $\L\neq 0$, hence
in this case, $M(\L,\succ)$ is irreducible.

On the other hand, if $\L=0,$ then it is clear that $$
M'(0,0)=\sum\limits_{k>0,\, \a_{1},\cdots,\a_{k}\in G_{+} } \F
L_{-a_1,i_1}\cdots L_{-\a_k,i_k}v_0 $$ is a proper
$U(\BB(G))$-submodule. Assume that for all $x,y\in G_+$ there exists
a positive integer $n$ such that $nx\succcurlyeq y$. By Steps 1--4,
  there exists $\varepsilon'\in G_+$ such that
$L_{-n\varepsilon',-1}v_0\in M'(0,0)\mbox { for all } n\in \N.$ Thus
for any $z\in G_+$, using $y=n\varepsilon' \succcurlyeq z$ for some
$n\in\N$, we
have
$$L_{-z,k-1}v_{0}=-((k+1)y)^{-1}L_{y-z,k}L_{-y,-1}v_{0}\in
M'(0,0)\mbox{ \ \ for all \ }k\in \Z_+.$$ We see that $M'(0,0)$ is
in fact an irreducible $\BB(G)$-module.

If there exists $x,y\in G_+$ such that $\N x\prec y$, then
$B(x)\prec y$ (cf.~(\ref{denote-Ba})). It is easy to verify that
$$W'=U(\BB(G))\{L_{-z,k}\,\,|\,\,z\in B(x),k\geq -1 \}v_\L,$$ is a
proper submodule of $M'(0,0)$ since $L_{-y}v_\L\notin W'$. \vskip
5pt
 (2) Suppose the order `` $\succ$"  is discrete (recall
 (\ref{discrete})).
 Note that $a\Z\subseteq G$. For
any $x\in G$, we write $x\succ a\Z$ if $x\succ na$ for any $n\in\Z$.
Let $$H_{+}=\bigl\{x\in G\bigm|x\succ a\Z\bigr\},\ \ \
H_{-}=-H_{+}.$$ Denote by $\BB(H_+)$ the subalgebra of $\BB(G)$
generated by $\{L_{\a,k}\,|\,\a\in H_+,k\geq -1\}.$ It is not
difficult to see that $G=a\Z\cup H_{+}\cup H_{-}$. Obviously,
$\BB(H_+)M_{a}(\L,\succ)=0$. Since $$ M(\L,\succ)\cong
U(\BB(G))\otimes_{U(\BB(a\Z)+\BB(H_+))}M_{a}(\L,\succ),
$$ it follows that the irreducibility of $\BB(G)$-module
$M(\L,\succ)$ imply the irreducibility of $\BB(a\Z)$-module
$M_{a}(\L,\succ)$.

Conversely, suppose $M_{a}(\L,\succ)$ is an irreducible
$\BB(a\Z)$-module. Let $u_0\notin \F v_\L$ be any given weight
vector in $M(\L,\succ)$. Then $u_0\in V_{r}\setminus{V_{r-1}}$ for
some $r\in\Bbb{N}$. We want to prove that $U(\BB(G))u_{0}\cap
M_{a}(\L,\succ)\neq\{0\}$, from which the irreducibility of
$M(\L,\succ)$ as a $\BB(G)$-module follows immediately.

Write $$u_0\equiv \sum\limits_{ {\stackrel{\sc
\a_{1}',\cdots,\a_{s}'\in H_{+},\, \a_{1},\cdots,\a_{r-s}\in
a\Z_{+},}{ \a_{1}'\succcurlyeq\cdots\succcurlyeq \a_s',\,
\a_{1}\succcurlyeq\cdots\succcurlyeq \a_{r-s}}}}
a_{\bar{\a},\bar{j}}L_{-\a_1',j'_1}\cdots
L_{-\a'_s,j'_s}L_{-\a_1,j_1}\cdots L_{-\a_{r-s},j_{r-s}}v_\L\mbox{
}({\rm mod}\mbox{ }V_{r-1}),$$ for some $a_{\bar{\a},\bar{j}}\in\F$,
where $j_t\geq j_{t+1}$ if $\a_t=\a_{t+1}$, and $j'_t\geq j'_{t+1}$
if $\a'_t=\a'_{t+1}$, and
$$(\bar{\a},\bar{j})=(\a'_1,\cdots,\a'_s,\a_1,\cdots,\a_{r-s},j_{1}',\cdots,j_{s}',j_{1},\cdots,
j_{r-s}).$$  Let
$$\bar{I}=\bigl\{(\a'_1,\cdots,\a'_s,\a_1,\cdots,\a_{r-s},j_{1}',\cdots,j_{s}',j_{1},\cdots,
j_{r-s}) \,|\,a_{\bar{\a},\bar{j}}\neq 0\bigr\}.$$ By assumption,
$\bar{I}\neq\emptyset$.

For any
\begin{eqnarray*}
(\bar{\a},\bar{j})&=&(\a'_1,\cdots,\a'_s,\a_1,\cdots,\a_{r-s},j_{1}',\cdots,j_{s}',j_{1},\cdots,
j_{r-s})\in \bar{I},\\
(\bar{\gamma},\bar{l})&=&(\gamma'_1,\cdots,\gamma'_t,\gamma_1,\cdots,\gamma_{r-t},l_{1}',\cdots,l_{t}',l_{1},\cdots,
l_{r-t})\in \bar{I},
\end{eqnarray*}we define $(\bar{\a},\bar{j})\succ'
(\bar{\gamma},\bar{l})$ if and only if
(cf.~(\ref{order1})--(\ref{order1+}))
\begin{eqnarray*}
(\a_{r-s},{\sc\cdots},\a_1,\a'_s,{\sc\cdots},\a'_1,j_{r-s},{\sc\cdots},j_{1},j_{s}',{\sc\cdots},
j_{1}')\succ(\gamma_{r-t},{\sc\cdots},\gamma_1,\gamma'_t,{\sc\cdots},\gamma'_1,l_{r-t},{\sc\cdots},l_{1},l_{t}',{\sc\cdots},
l_{1}').
\end{eqnarray*}
 Let
$$
(\bar{\b},\bar{i})=(\b'_1,\cdots,\b'_1,\b'_{t+1},\cdots,\b'_m,\b_1,\cdots,\b_{r-m},i_1',\cdots,i_t',i_{t+1}',\cdots,
i_{m}',i_{1},\cdots,i_{r-m})
$$
with $\b'_1\neq \b'_{t+1}$, be the unique maximal element in
$\bar{I}$ with respect to $\succ'$. Note that
$\b'_{1}-\a'_{k}-a\succcurlyeq 0$ if $\b'_{1}\neq \a'_{k}.$ Then for
$k_1\gg 0,$ we have
\begin{eqnarray*}
u(1):=L_{\b_1'-a,k_1}u_{0}\equiv\!\!\!\!&
\mbox{$\sum\limits_{{\stackrel{\sc \a_{1}',\cdots,\a_{s-1}'\in
H_{+},\, \a_{1},\cdots,\a_{r-s}\in a\Z_{+},}{
\a_{1}'\succcurlyeq\cdots\succcurlyeq \a_{s-1}',\,
\a_{1}\succcurlyeq\cdots\succcurlyeq \a_{r-s},\, k'_1\in\N}}}$}
a_{\bar{\a},\bar{j}}^{(1)}L_{-\a_1',j_1'}\cdots
L_{-\a'_{s-1},j_{s-1}'}
\\\!\!\!\!\!\!\!\!\!\!\!\!\!\!\!\!&\times\  L_{-a,k_1'}L_{-\a_1,j_1}\cdots
L_{-\a_{r-s},j_{r-s}}v_\L\ \ \ ({\rm mod}\,{V_{r-1}}),
\end{eqnarray*} for some $a_{\bar{\a},\bar{j}}^{(1)}\in\F$.
Set
$$
\begin{array}{l} \bar{I}^{(1)}=\{(\a'_1,\cdots,\a'_{s-1},a,\a_1,\cdots,\a_{r-s},j_{1}',\cdots,j_{s-1}',
k_1',j_1,\cdots,j_{r-s})\,|\,a_{\bar{\a},\bar{j}}^{(1)}
\neq 0\} .\\[9pt]
(\bar{\b},\bar{i})^{(1)}\!=\!(\b'_1,{\sc\cdots},\b'_{1},\b'_{t+1},
{\sc\cdots},\b'_{m},a,\b_1,{\sc\cdots},\b_{r-m},i_1',{\sc\cdots},i_{t-1}',
i_{t+1}',{\sc\cdots},i_{m}',k_1+i'_t,i_{1},
{\sc\cdots},i_{r-m}).\end{array}$$
The term $$L_{-\b'_1,i_1'}\cdots
L_{-\b'_{1},i_{t-1}'}L_{-\b'_{t+1},i_{t+1}'}\cdots
L_{-\b'_{m},i_{m}'}L_{a,k_1+i'_t}L_{\b_1,i_{1}} \cdots
L_{-\b_{r-m},i_{r-m}}v_\L$$ appears in $u(1)$ since the
corresponding coefficient is
$$-m((k_1+1)\b_1'+(i_s'+1)(\b_1'-a))a_{\bar{\b},\bar{i}}\neq 0\mbox{ \ \ for
some \ }m\in \N.$$ Thus $\bar{I}^{(1)}\neq \emptyset.$ Now for
$l=2,\cdots,r$, we define recursively and prove by induction that

(i) Choose $k_l\gg 0$ and let $u(l)=L_{\b_{m-l+1}'-a,k_l}u(l-1)$.
Then
\begin{eqnarray*}u(s)&\equiv\!\!\!\!&\sum\limits_{{\stackrel{\sc
k_1',\cdot,k_l'\in \N,\, \a_{1}',\cdots,\a_{s-l}'\in H_{+},}
{\a_{1},\cdots,\a_{r-s}\in a\Z_{+}}}}
a_{\bar{\a},\bar{j}}^{(l)}L_{-\a_1',j'_1}\cdots
L_{-\a'_{s-l},j_{s-l}'}\\&&\times\ L_{-a,k_l'}\cdots
L_{-a,k_1'}L_{-\a_1,j_1}\cdots L_{-\a_{r-s},j_{r-s}}v_\L ~({\rm
mod}~{V_{r-1}}),
\end{eqnarray*}
for some $a_{\bar{\a},\bar{j}}^{(l)}\in\F$, where
$\a'_1\succcurlyeq\cdots \succcurlyeq\a'_{s-l},\,
\a_1\succcurlyeq\cdots\succcurlyeq \a_{r-s}.$

(ii) Let
$$\bar{I}^{(l)}=\{(\a_{1}',\cdots,\a'_{s-l},a,\cdots,a,\a_1,\cdots,\a_{r-s},j'_1,\cdots,j_{s-l}',k_l',\cdots,k_1',
j_{1},\cdots,j_{r-s})\,|\,a_{\bar{\a},\bar{j}}^{(s)}\neq 0\}.$$\\
Then $\bar{I}^{(l)}\neq \emptyset.$

Now letting $l=m$ and noting that $u(m)$ is a weight vector, we
obtain that $0\neq u(m)\in U(\BB(G))u_{0}\cap M_{a}(\L,\succ)$ as
required.\hfill$\Box$\vskip7pt

\cl{\bf\S3. \ Verma modules over $\BB(\Z)$
}\setcounter{section}{3}\setcounter{theo}{0}\setcounter{equation}{0}

Following [S1], we realize the Lie algebra $\BB(\Z)$ in the space
$\F[x,x^{-1},t]\oplus \F c$ with the bracket
\begin{equation}\label{e-2.1}
[x^\a f(t),x^\b g(t)]=x^{\a+\b}(\b f'(t)g(t)-\a
f(t)g'(t))+\a\d_{\a,-\b}f(0)g(0)c,
\end{equation}for $\a,\b\in \Z,f(t),g(t)\in \F[t],$ where the
prime stands for the derivative $\frac{d}{d{\ssc\,}t}.$

We denote $$L_{\a,i}=x^\a t^{i+1} \mbox{ for }\a\in\Z, i\geq -1.$$
Then (\ref{e-2.1}) is equivalent to (\ref{def}).

We always use the normal order on $\Z$. Denote by $M(\L)$ the Verma
$\BB(\Z)$-module  with highest weight vector $v_\L.$ Suppose $M(\L)$
is reducible. Let $M'$ denote the maximal proper submodule of
$M(\L),$ and set $L(\L)=M(\L)/M'$, the irreducible highest weight
module of weight $\L$. Set
$$\AA=\{a\in\BB(\Z)\,|\,av_\L\in M'\} \mbox{ \ \ and \ \
}\PP=\AA+\BB(\Z)_0.$$ Clearly, $\BB(\Z)_+\subset \AA,$ and $\PP$ is
a subalgebra of $\BB(\Z).$

\begin{lemm}\label{lemm}
$\PP$ is a parabolic subalgebra of $\BB(\Z)$, namely,
\begin{equation}\label{equa-lemm}\PP\supset
\BB(\Z)_0+\BB(\Z)_+\neq\PP.\end{equation}
\end{lemm}
\noindent{\it Proof.~} The proof of (\ref{equa-lemm}) is equivalent
to proving
\begin{equation}\label{equa-lemm+}
\mbox{ $\PP\cap\BB(\Z)_m\ne 0$ \ \ for some \ \ $m<0$.
}\end{equation}
 Let $n$ be the
minimal positive integer such that $U(\BB(\Z))_{-n}v_\L\cap M'\ne
0.$
If $n=1,2$, one can easily verify that (\ref{equa-lemm+}) holds.
Assume that $n>2.$ Then there exists a vector $u$ of weight $-n$ in
$M'.$ Write
$$ u=\sum\limits_{{\stackrel{^{^{\sc 1\leq l\leq n,\,\a_1\leq\cdots\leq
\a_l}}}{_{_{\sc \a_1+\cdots
\a_l=-n}}}}}c_{\ul{\a},\ul{j}}x^{-\a_1}t^{j_1}\dots
x^{-\a_l}t^{j_l}v_\L\in M'\mbox{ \ \ for some \
}c_{\ul{\a},\ul{j}}\in\F,$$ where
$\ul{\a}=(\a_1,\cdots,\a_l),\,\ul{j}=(j_1,\cdots,j_l),$ and $j_s\leq
j_{s+1}$ if $\a_s=\a_{s+1}$ for $1\le s\le l-1$. Moreover, we denote
$$(\ul{\a},\ul{j})=(\a_1,\cdots,\a_l,j_1,\cdots,j_l),\mbox{ \ \ and \
\ }\ul{1}=(-1,\cdots,-1)\mbox{ \ ($n$ copies of $-1$'s)}.$$

\begin{clai} $c_{\ul{1},\ul{j}}\neq 0$ for some $\ul{j}.$
\end{clai}

Write (recall (\ref{gradation}))$$ u\equiv
\sum\limits_{{\stackrel{0<\a_1\leq\cdots\leq \a_r,}{\a_1+\cdots
\a_l=-n}}}c_{\ul{\a},\ul{j}}x^{-\a_1}t^{j_1}\dots
x^{-\a_r}t^{j_r}v_\L\mbox{ }(\rm mod ~V_{r-1}),$$where
$\ul{J}=\{(\ul{\a},\ul{j})\,|~c_{\ul{\a},\ul{j}}\neq 0\}\neq
\emptyset.$ Assume that there exists $(\ul{\a},\ul{j})\in \ul{J}$
such that $\ul{\a}\neq \ul{1}.$
Let$$(\ul{\b},\ul{i})=(1,\cdots,1,\b_s,\cdots,\b_r,i_1,\cdots,i_{s-1},i_s,\cdots,i_r)$$
be the unique maximal element in $\ul{J}$ (here we use the order
defined as in (\ref{order1+})), where $s\geq 1$ and $\b_s\neq 1.$ By
assumption, we have $s\neq r+1.$ Then for $k\gg 0,$
$$xt^k\cdot u\equiv
\sum\limits_{\a'_1\leq\cdots\leq
\a'_r}c_{\ul{\a}',\ul{j}'}x^{-\a'_1}t^{j'_1}\dots
x^{-\a'_r}t^{j'_r}v_\L\mbox{ }(\rm mod ~V_{r-1}).$$Set
$$\begin{array}{l}
\ul{J}'=\{(\ul{\a}',\ul{j}')=(\a'_1,\cdots,\a'_r,j'_1,\cdots,j'_r)\,|~c_{\ul{\a}',\ul{j}'}\neq
0\},\\[9pt](\ul{\b},\ul{j})'=(1,\cdots,1,\b_s-1,\b_{s+1},\cdots,\b_r,i_1,\cdots,i_{s-1},
k+i_s-1,i_{s+1},\cdots,i_r)\end{array}$$ such that
$(\ul{\b},\ul{j})'$ is the unique maximal element in $\ul{J}'$. The
term $$x^{-1}t^{i_1}\dots x^{-1}t^{i_{s-1}}
x^{-\b_s+1}t^{k+i_{s}-1}x^{-\b_{s+1}}t^{i_{s+1}}\cdots
x^{-\b_{r}}t^{i_{r}}v_\L$$ appears in  $xt^k\cdot u$ since the
corresponding coefficient  is $-m(\b_sk-i_s)c_{\ul{\b},\ul{i}}\neq
0$ for some $m\in \N.$ Thus $\ul{J}'\neq \emptyset$ and $0\neq
U(\BB(\Z))_{-n+1}v_\L\cap M',$ a contradiction with the assumption.
Our claim follows. \vs{6pt}

 Now we can write $$u\equiv \sum\limits_{i_1\leq\cdots\leq
i_n}c_{\ul{i}}x^{-1}t^{i_1}\dots x^{-1}t^{i_n}v_\L+
\sum\limits_{l_1\leq\cdots\leq l_{n-1}}c'_{\ul{l}}x^{-1}t^{l_1}\dots
x^{-1}t^{l_{n-2}}x^{-2}t^{l_{n-1}}v_\L\mbox{ }(\rm mod ~V_{n-2}),
$$ for some $c_{\ul i},c'_{\ul l}\in\F$, where $$\ul{I}'=\{\ul{i}=(i_1,\cdots,i_n)\,|~c_{\ul{i}}\neq 0\}\neq
\emptyset.$$ For any $\ul{i},\,\ul{i}'\in \ul{I}',$ we define
$\ul{i}> \ul{i}'$ as in (\ref{order-i}). Let
$\ul{j}=(j_1,\cdots,j_n)$ be the unique maximal element in
$\ul{I}'.$ For $k\gg 0,$ we have a nonzero weight vector
$$xt^k\cdot u=\sum\limits_{i'_1\leq\cdots\leq i'_{n-1}}d_{\ul{i}'}x^{-1}t^{i'_1}\cdots x^{-1}t^{i'_{n-1}}v_\L
\mbox{ }({\rm mod}\mbox{ }V_{n-2})\mbox{ \ \ for some \ }d_{\ul
i'}\in\F,$$ since the coefficient corresponding to
$$(1,\cdots,1,j_1,\cdots,j_{n-2},k+j_{n-1}+j_n-1)$$ is
$$m(k+j_{n-1})(k+j_{n-1}+j_n-1)c_{\ul{j}}-(2k+j_{n-1}+j_n)d_{\ul{h}}\neq
0\mbox{ \ \ for some \ }m\in\N,$$ where
$\ul{h}=(j_1,\cdots,j_{n-2},j_{n-1}+j_n).$ Thus $0\neq
U(\BB(\Z))_{-n+1}v_\L\cap M',$ a contradiction. Our lemma
follows.\hfill$\Box$\vskip5pt

By the lemma \ref{lemm}, we have $\PP_{-1}=\PP\cap\BB(\Z)_{-1}\neq
0.$ Let $f(t)$ be the monic ploynomial with minimal degree such that
$x^{-1}f(t)\in \PP.$ We shall call such polynomial {\it charactic
polynomial} (cf.~[S1, KL]). Set $a=x^{-1}f(t)$. Since $M'$ is a
proper submodule, we have $b\cdot av_\L=0$ for any $b\in \BB(\Z)_+.$
From (\ref{e-2.1}), we have
$$[xg(t),x^{-1}f(t)]v_\L=\L(f'(t)g(t)+f(t)g'(t)-f(0)g(0)c)=0~\mbox{
\ for all \ }g(t)\in \F[t].$$

A weight $\L\in \BB(\Z)_0^*$ is described by the {\it central
charge} $c=\L(c)$ and its {\it label} $\L_i=\L(t^i)$ for $i\geq 0.$
We introduce the {\it generating series}
$$\Delta_\L(z)=c+\sum\limits_{i=0}^\infty
{\frac {z^{i+1}}{z!}}\L_i=c-\L(ze^{zt}).$$ From [S1] and the above
arguments, we obtain the following theorem.
\begin{theo}\label{theo2} The following conditions are equivalent:
\begin{itemize}\parskip-6pt

\item[{\rm(1)}] $M(\L)$ is reducible.

\item[{\rm(2)}] $\PP_{-1}\neq \{0\}.$

\item[{\rm(3)}] $\Delta_\L(z)$ is a {\it quasipolynomial}.

\item[{\rm(4)}] $L(\L)$ is  {\it quasifinite}.
\end{itemize}
\end{theo}

Now Theorem \ref{theo1}(1) follows from Theorem \ref{theo2}.
 \vskip12pt

\cl{\bf References}\vs{0pt}

\vskip5pt\small
\parindent=8ex\parskip=2pt\baselineskip=2pt

\re{B} R. Block, On torsion-free abelian groups and Lie algebras,
  {\it Proc. Amer. Math. Soc.} {\bf 9} (1958), 613--620.

\re{DZ} D. Dokovic, K. Zhao, Derivations, isomorphisms and
cohomology of generalized Block algebras, {\it Algebra Colloq.} {\bf
3} (1996), 245--272.


\re{KL} V. Kac, J. Liberati, Unitary quasi-finite representations of
$W_\infty$, {\it Lett. Math. Phys.} {\bf 53} (2000), 11--27.

\re{KR} V. Kac, A. Radul, Quasi-finite highest weight modules over
the Lie algebra of differential operators on the circle,{\it Comm.
Math. phys.} {\bf 157} (1993), 429--457.

\re{LT}  W. Lin, S.  Tan, Nonzero level Harish-Chandra modules over
the Virasoro-like algebra, {\it J. Pure Appl. Algebra}, in press.




\re{S1} Y. Su, Quasifinite representations of a Lie algebra of Block
type, {\it J. Algebra} {\bf 276} (2004), 117--128.

\re{S2} Y. Su, Quasifinite representations of a family of Lie
algebras of Block type, {\it J. Pure Appl. Algebra} {\bf192} (2004),
293--305.


\re{WZ} X. Wang, K. Zhao, Verma modules over the Virasoro-like
algebra, {\it J. Australia Math.}, in press.

\re{X1} X. Xu, Generalizations of Block algebras, {\it Manuscripta
Math}. {\bf 100} (1999), 489--518.

\re{X2} X. Xu, Quadratic conformal superalgebras, {\it J. Algebra}
 {\bf 231} (2000), 1--38.


\re{ZM} L. Zhu, D. Meng, Structure of degenerate Block algebras,
{\it Algebra Colloq.} {\bf 10} (2003), 53--62.


\end{document}